\documentclass[11pt]{article}
\usepackage{mathrsfs}
\usepackage{amssymb}
\usepackage{amsfonts}
\usepackage{color,xcolor}
\usepackage{graphicx}
\usepackage{epstopdf}
\usepackage{geometry}
\usepackage{manfnt}
\usepackage[pagewise]{lineno}
\RequirePackage[colorlinks,citecolor=blue,urlcolor=blue]{hyperref}
\newtheorem{theorem}{Theorem}[section]

\newtheorem{lemma}{Lemma}[section]

\newtheorem{remark}{Remark}[section]

\def\[{{\Big[}}\def\]{{\Big]}}\def\({{\Big(}}\def\){{\Big)}}
\def\={&\!\!=\!\!&}
\def\geq{\geqslant}\def\leq{\leqslant}

\begin{document}
\title{\bf Notes on spatial twisted central configurations for $2N$-body problem}
\author{Liang Ding$^{1}$, Juan Manuel S\'{a}nchez-Cerritos$^{2,3}$ and Jinlong Wei$^{4\,\ast}$
\\ {\small \it $^1$School of Data Science and Information Engineering, Guizhou} \\  {\small \it Minzu University, Guiyang, 550025, China}\\
{\small \tt ding2016liang@126.com}\\  {\small \it $^2$College of Mathematics and Statistics, Chongqing Technology} \\  {\small \it and Business University, Chongqing 400067, China}\\
 {\small \it
$^3$Chongqing Key Laboratory of Social Economy and Applied Statistics} \\  {\small \it Chongqing 400067, China}\\
 {\small \tt sanchezj011@outlook.com}\\  {\small \it $^4$School of Statistics and Mathematics, Zhongnan University} \\  {\small \it of Economics and Law, Wuhan, 430073, China}\\
 {\small \tt weijinlong.hust@gmail.com}}
\date{}
 \maketitle
\noindent{\hrulefill} \vskip1mm\noindent
 {\bf Abstract}
We study the spatial central configuration formed by two twisted regular $N$-polygons. For any twist angle $\theta$ and any ratio of the masses $b$ in the two regular $N$-polygons,  we prove that the sizes of the two regular $N$-polygons must be equal.

\vskip2mm\noindent {\bf Keywords.}  Spatial twisted $2N$-body problem; Central configurations;
The radio of sizes; The radio of masses; Twisted angle

\vskip2mm\noindent {\bf MSC (2010):} 70F07, 70F15

 \vskip1mm\noindent{\hrulefill}
\section{Introduction} \label{sec1}
\setcounter{equation}{0}
The study of central configurations is a very important subject in celestial mechanics with a long and varied
history \cite{Moeckel1990}, and a well-known result is that finding the relative equilibrium solutions of the classical $N$-body problem and the planar central configurations is equivalent \cite{Hampton2006}. The numbers and shapes of central configurations for the Newtonian $N$-body $(N\geq4)$ problem are important and difficult problems in celestial mechanics \cite{Saari1980}. In  \cite{Smale}, Smale took it as one of the most important 18 mathematical problems (the sixth one) for the 21st century.
Though there are a lot of elegant works on central configurations \cite{Albouy2008,Albouy2012,Fernandes2017,Hampton2005,Perko1985,
Wang2019,Yu2012,Zhang2003,Zhang2018}, it is a few works to find concrete central configurations since it is a difficult problem. Firstly, we give some preliminaries on spatial central configuration.

Given $2N$ mass points $m_k$ ($k=1,\ldots,2N$) with position $q_k \in \mathbb{R}^{3}$. Denote $r_{kj}=|q_k-q_j|$ as the Euclidean distance between the mass particles $m_k$ and $m_j$ and let
$\mathbf{q}=(q_1,\ldots,q_{2N})$ $\in \mathbb{R}^{6N}$, the center of mass of the system is $c_{0}=\frac{1}{M}\sum_{k=1}^{2N}m_kq_k$,
where $M=m_{1}+\ldots +m_{2N}$ is the total mass.
In the inertial system, the motion equations of $2N$ bodies can be described by Newton's three laws on classical mechanics and Newton's universal gravitation law:
$$
m_k\ddot{q_k}=\frac{\partial (\sum_{1\leq s<j\leq 2N}\frac{m_{j}m_s}{r_{js}})}{\partial q_k },\ \ \ k=1,\ldots,2N.
$$

For the above Newtonian $2N$-body problem with configuration $\mathbf{q} \in \mathbb{R}^{6N}$, we now give the following definition.

\medskip \noindent
\textbf{Definition 1.1.}  (\cite{Wintner1947})
Given $2N$ mass points $m_k$ with position $q_k \in \mathbb{R}^{3}$,
$k=1,\ldots, 2N$. A configuration $\mathbf{q}=(q_1,\ldots,q_{2N})^{T}\in (\mathbb{R}^{3})^{2N}$ is called a central configuration if $q_{k}\neq q_{j} \,\, when \,\, k\neq j$, and there exists a constant $\lambda\in \mathbb{R}$ such that
\begin{eqnarray}\label{Equ1.1}
\left\{\begin{array}{ll}
\sum\limits_{{j\neq k\atop{{1\leq j\leq 2N}}}}
\frac{m_j m_k}{|q_j-q_k|^{3}}(q_j-q_k)=-\lambda m_{k}(q_{k}-c_{0}), \quad \  \  k=1,\ldots,2N, \\
\ \lambda=\frac{U(\mathbf{q})}{I(\mathbf{q})},
\end{array}\right.
\end{eqnarray}
where the Newtonian potential $V$ is given by
\begin{eqnarray*}
V(\mathbf{q})=-U(\mathbf{q})=-\sum_{1\leq k<j\leq 2N}\frac{m_j m_k}{|q_j-q_k|},
\end{eqnarray*}
and the moment of inertia of $\mathbf{q}$ is given by
\begin{eqnarray*}
 I(\mathbf{q})=\sum_{k=1}^{2N}m_k|\mathbf{q_k}-c_{0}|^{2}.
\end{eqnarray*}

In this paper, we study the two twisted regular $N$-polygonal central configurations, and we use the following notation: suppose two parallel regular polygons, one regular $N$-polygon and the other regular $N$-polygon with
distance $h>0$ are placed in $\mathbb{R}^{3}$ (See Figure 1). Assume that the particles $q_1,\ldots,q_N$ with the same mass $m$ locate at the vertexes of one regular $N$-polygon; the particles $q_{N+1},\ldots,q_{2N}$ with the same mass $bm$ locate at the vertexes of the other regular $N$-polygon.
Let $\rho_{k}=e^{i\theta_{k}}$ be the $k$-th root of the $N$-roots of unity,
where $\theta_{k}=2k\pi/N (k=1,\ldots, N)$, and $\rho_{l}=a\rho_{k}\cdot e^{i\theta}$ ($l=k+N$),
where $a>0$, and $0\leq\theta\leq2\pi$ is called the twisted angle. So the ratio of the masses is $b$, and the ratio of the sizes of the
two polygons is $a$. Moreover, let the coordinates of the particles $q_{1},\ldots,q_{N}$ and
$q_{N+1},\ldots, q_{2N}$ be $q_{k}=(\rho_{k}, 0)$ $(k=1,\ldots, N)$, $q_{l}=(\rho_{l}, h)$ $ (l=N+1,\ldots, 2N)$, respectively.

Suppose the configuration is formed by two twisted regular $N$-polygons ($N\geq2$) with distance $h>0$. Moeckel and Simo obtained that: if $N<473$ and the twist angle $\theta=0$, there is a unique pair of spatial central configurations of regular $N$-polygons \cite[Thereom 2]{Moeckel1995}, and if $N\geq473$ and the twist angle $\theta=0$, there exists no spatial central configuration for $b<\mu_{0}(N)<1$ \cite[Thereom 2, and Proposition 3 and Page 986, lines 1-2]{Moeckel1995}. This result is generalized by Zhang and Zhu \cite{Zhang2002} to $\theta=\pi/N$ for the two twisted regular $N$-polygons ($N\geq2$) with distance $h>0$. For every $N\geq 2$, they proved that if $b=1$, $a=1$ $\theta=\pi/N$,  and the configuration formed by two twisted regular $N-$polygons with distance $h>0$ is a central configuration, then there exists only one $h>0$ such that $q_{1},\ldots,q_{N},q_{N+1},\ldots,q_{2N}$ form a spatial central configuration. For more details in these direction, one can refer to \cite{Wang2015,Xie2000,Xie2012,Yu2015,Zhang2003}.

Note that in 2003, for the planar twisted central configurations (i.e. $h=0$) formed by two regular $N-$polygons with any twist angle $\theta$, Zhang and Zhou \cite{Zhang2003} arrived at the conclusion that the values of masses in each separate regular $N$-polygons must be equal (but without detailed proof); in 2015, based on the eigenvalues of circulant matrices, Wang and Li \cite{Wang2015} investigated the masses of the $2N$ bodies for $h\geq0$ with twist angle $\theta=0$, and
they also obtained the values of masses in each separate regular $N$-polygons must be equal. Moreover, we also note that Yu and Zhang \cite{Yu2012} proved that if the central configuration is formed by two twisted regular $N$-polygons with distance $h\geq0$, then the twist angles must be $\theta=0$ or $\theta=\pi/N$, so
 we want to study that for the spatial central configuration ($h>0$) formed by two twisted regular $N$-polygons with twist angle $\theta=\pi/N$, whether the sizes of each separate regular $N$-polygons must be equal or not ?

 In this paper, we attempt to answer this question.

\section{Main result} \label{sec2}
\setcounter{equation}{0}
By analysing the relationship between the ratio of the masses of the two regular $N$-polygons $b$ and the ratio of the sizes $a$, we obtain the following main result:

\begin{theorem} \label{the2.1} For the spatial twist central configuration formed by two twisted regular $N$-polygons ($N\geq3$) with any twist angle $\theta$ and any ratio of the masses $b$, then the sizes of the two regular $N$-polygons must be equal.
\end{theorem}
\begin{figure}
\centering
\scalebox{0.5}[0.5]{
\includegraphics{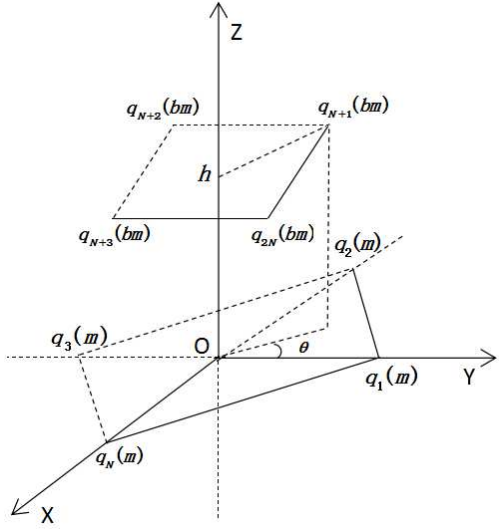}}
 \caption{Theorem 2.1}
\label{figure 1}
\end{figure}

\begin{remark}\label{rem2.1}
From the above theorem, we know for any ratio of masses $b\in\mathbb{R}$, we have $a=1$. Then by the following Remark \ref{rem3.1}, it arrives at an interesting result that $b=1$, which means for the spatial central configuration formed by two twisted regular $N$-polygons with any twist angle $\theta$,
the masses of the $2N$ bodies must be equal to each other.
\end{remark}

\section{Some useful lemmas}\setcounter{equation}{0}

Set
\begin{eqnarray}\label{3.1}
\left\{
  \begin{array}{ll}
 x=\sum_{1\leq k\leq N-1}\frac{1-\rho_{k}}{|1-\rho_{k}|^{3}},
   \\ [3mm]
y=\sum_{k=1}^{N}\frac{\cos(\theta_{k}+\theta)}{[1+a^{2}-2a\cos(\theta_{k}+
\theta)+h^{2}]^{\frac{3}{2}}},\\ [3mm]
z=\sum_{k=1}^{N}\frac{1}{[1+a^{2}-2a\cos(\theta_{k}+\theta)+h^{2}]^{\frac{3}{2}}},
  \end{array}
\right.
\end{eqnarray}
and before proving the main result, we introduce some lemmas which will serve us well later.
\begin{lemma} \cite [Theorem 1.8]{Yu2012} \label{lem3.1} If the central configuration is formed by two twisted regular $N$-polygons ($N\geq2$) with distance $h\geq0$,
then the twist angle only $\theta=0$ or $\theta=\pi/N$.
\end{lemma}

\begin{lemma}  \cite[Corollary 1.10]{Yu2012} \label{lem3.2}
The spatial configuration formed by two twisted regular $N$-polygons ($N\geq2$)
with distance $h>0$, is a central configuration if and only if the parameters $a$, $b$ and $h$ satisfy the following relationships:
\begin{eqnarray}\label{3.2}
&&bay=x-z
\end{eqnarray}
and
\begin{eqnarray}\label{3.3}
\frac{b}{a^{2}}x-baz=y.
\end{eqnarray}
\end{lemma}

\begin{remark}\label{rem3.1}
From Lemma \ref{lem3.2}, we see that if the configuration formed by two twisted regular $N$-polygons ($N\geq2$)
with the ratio of size $a=1$ and distance $h>0$, is a central configuration, then $b=1$. But from $b=1$ and $h>0$, we can not obtain $a=1$ directly.
\end{remark}

\begin{lemma} \label{lem3.3} For any distance $h>0$, whether the twist angle $\theta=0$ or $\theta=\pi/N$, we have the following equality
\begin{eqnarray}\label{3.4}
\sum_{k=1}^{N}\frac{\cos(\theta_{k}-\theta)}{[1+a^{2}-2a\cos(\theta_{k}-\theta)+h^{2}]^{\frac{3}{2}}}=
\sum_{k=1}^{N}\frac{\cos(\theta_{k}+\theta)}{[1+a^{2}-2a\cos(\theta_{k}+\theta)+h^{2}]^{\frac{3}{2}}}.
\end{eqnarray}
\end{lemma}
\textbf{Proof.} For the case of $\theta=0$, equality (\ref{3.4}) obviously holds. For the rest case of $\theta=\pi/N$, by $\theta_{k}=2k\pi/N (k=1,\ldots, N)$, we have
\begin{eqnarray*}
&&\sum_{k=1}^{N}\frac{\cos(\theta_{k}-\theta)}{[1+a^{2}-2a\cos(\theta_{k}-\theta)+h^{2}]^{\frac{3}{2}}}=
\sum_{k=1}^{N}\frac{\cos\frac{(2k-1)\pi}{N}}{[1+a^{2}-2a\cos\frac{(2k-1)\pi}{N}+h^{2}]^{\frac{3}{2}}}\cr
&=&\sum_{k=1}^{N}\frac{\cos\frac{(2k+1)\pi}{N}}{[1+a^{2}-2a\cos\frac{(2k+1)\pi}{N}+h^{2}]^{\frac{3}{2}}}
=\sum_{k=1}^{N}\frac{\cos(\theta_{k}+\theta)}{[1+a^{2}-2a\cos(\theta_{k}+\theta)+h^{2}]^{\frac{3}{2}}},
\end{eqnarray*}
which implies that whether $\theta=0$ or $\theta=\pi/N$, equality (\ref{3.4}) holds. $\Box$

\begin{lemma} \cite[Theorem 1.2]{Zhang2002} \label{lem3.4}
The configuration formed by two twisted regular $N$-polygons ($N\geq2$) with distance $h>0$ is a central configuration if and only if the parameters $a$, $b$, $h$, $\theta$ satisfy the following relationships:
\begin{eqnarray}\label{3.5}
\lambda\frac{N}{M}=\frac{1}{1+b}\(x+\sum_{k=1}^{N}\frac{b(1-ae^{i\theta}\rho_{k})}{\[|1-ae^{i\theta}\rho_{k}|^{2}+h^{2}\]^{\frac{3}{2}}}\),
\end{eqnarray}
and
\begin{eqnarray}\label{3.6}
\lambda\frac{N}{M}=\sum_{k=1}^{N}\frac{1}{\[|1-ae^{i\theta}\rho_{k}|^{2}+h^{2}\]^{\frac{3}{2}}},
\end{eqnarray}
and
\begin{eqnarray}\label{3.7}
\lambda\frac{N}{M}=\frac{e^{-i\theta}}{a(1+b)}\(\sum_{1\leq k\leq N-1}\frac{b(1-\rho_{k})e^{i\theta}}{a^{2}|1-\rho_{k}|^{3}}+\sum_{k=1}^{N}\frac{ae^{i\theta}-\rho_{k}}{\[|ae^{i\theta}-\rho_{k}|^{2}+h^{2}\]^{\frac{3}{2}}}\),
\end{eqnarray}
where $M=m_{1}+\ldots +m_{2N}$.
\end{lemma}

\section{Proof of Theorem \ref{the2.1}}\setcounter{equation}{0}
From Lemma \ref{lem3.1}, the twist angle must be $\theta=0$ or $\theta=\pi/N$, and in the following, we use $\theta$ to represent the twist angle $0$ or $\pi/N$.

By the symmetry of the configuration, in the following steps, we need only consider the ratio of masses $0<b\leq1$, and we divide the proof into six steps.

\textbf{Step 1.} Let $h>0$,
we claim that $y=y(h)>0$.

In fact, denoting $\sum_{i=1}^{N}\rho_{k}=A$, it is easy to see that $A=0$.

Observing that
\begin{eqnarray*}
\sum^{N}_{k=1} e^{i\theta_k}=\sum^{N}_{k=1}\rho_k=A=0,
\end{eqnarray*}
therefore
\begin{eqnarray*}
e^{-i\theta}\sum^{N}_{k=1}\rho_k=\sum^{N}_{k=1} e^{i(\theta_k-\theta)}=0,
\end{eqnarray*}
and it means that
\begin{eqnarray}\label{4.1}
\sum^{N}_{k=1}\cos(\theta_k-\theta)=\sum_{ k\in I_1}\cos(\theta_k-\theta)+\sum_{ k\in I_2}\cos(\theta_k-\theta)=0,
\end{eqnarray}
where
\begin{eqnarray*}
I_1=\{1\leq k\leq N \ | \ \cos(\theta_k-\theta)>0 \},
\end{eqnarray*}
and
\begin{eqnarray*}
I_2=\{1\leq k\leq N \ | \ \cos(\theta_k-\theta)\leq0\}.
\end{eqnarray*}

Since $N\geq 3$, $I_1\neq \varnothing$. For any $k_1\in I_1$ and any $k_2\in I_2$, since $h>0$, then
\begin{eqnarray*}
\cos(\theta_{k_1}-\theta)> 0, \ -\cos(\theta_{k_2}-\theta)\geq 0,
\end{eqnarray*}
and
\begin{eqnarray}\label{4.2}
\frac{1}{[1+a^{2}-2a\cos(\theta_{k_1}-\theta)+h^{2}]^{\frac32}}>\frac{1}{[1+a^{2}-2a\cos(\theta_{k_2}-\theta)+h^{2}]^{\frac32}}.
\end{eqnarray}
By (\ref{4.1}), we have
\begin{eqnarray*}
\sum_{ k\in I_1}\cos(\theta_k-\theta)=-\sum_{ k\in I_2}\cos(\theta_k-\theta).
\end{eqnarray*}
With the aid of (\ref{4.2}),
we arrive at
\begin{eqnarray*}
&&\sum_{ k\in I_1}\frac{\cos(\theta_k-\theta)}{[1+a^{2}-2a\cos(\theta_{k}-\theta)+h^{2}]^{\frac32}}\cr
&\geq&
\[\sum_{ k\in I_1}\cos(\theta_k-\theta)\]\cdot\min_{k\in I_1}\frac{1}{[1+a^{2}-2a\cos(\theta_{k}-\theta)+h^{2}]^{\frac32}}\cr
&=&\(\sum_{ k\in I_2}[-\cos(\theta_k-\theta)]\)\cdot\min_{k\in I_1}\frac{1}{[1+a^{2}-2a\cos(\theta_{k}-\theta)+h^{2}]^{\frac32}}\cr
&>&\(\sum_{ k\in I_2}[-\cos(\theta_k-\theta)]\)\cdot\max_{k\in I_2}\frac{1}{[1+a^{2}-2a\cos(\theta_{k}-\theta)+h^{2}]^{\frac32}}\cr
&\geq&\sum_{ k\in I_2}\frac{-\cos(\theta_k-\theta)}{[1+a^{2}-2a\cos(\theta_{k}-\theta)+h^{2}]^{\frac32}}.
\end{eqnarray*}
Hence
\begin{eqnarray}\label{4.3}
&&\sum_{ k\in I_1}\frac{\cos(\theta_k-\theta)}{[1+a^{2}-2a\cos(\theta_{k}-\theta)+h^{2}]^{\frac32}}
+\sum_{ k\in I_2}\frac{\cos(\theta_k-\theta)}{[1+a^{2}-2a\cos(\theta_{k}-\theta)+h^{2}]^{\frac32}}\cr
&=&\sum^{N}_{k=1}\frac{\cos(\theta_{k}-\theta)}{[1+a^{2}-2a\cos(\theta_{k}-\theta)+h^{2}]^{\frac32}}>0.
\end{eqnarray}
Moreover, by (\ref{3.4}) in Lemma \ref{lem3.3}, then
\begin{eqnarray}\label{4.4}
&&y=\sum^{N}_{k=1}\frac{\cos(\theta_{k}+\theta)}{[1+a^{2}-2a\cos(\theta_{k}+\theta)+h^{2}]^{\frac32}}
\nonumber\\&=&\sum^{N}_{k=1}\frac{\cos(\theta_{k}-\theta)}{[1+a^{2}-2a\cos(\theta_{k}-\theta)+h^{2}]^{\frac32}}.
\end{eqnarray}
In view of (\ref{4.3}) and (\ref{4.4}), we conclude that $y(h)>0$ for $h>0$.

\textbf{Step 2.} We prove that for the spatial twist $2N$-body problem, if the ratio of the masses $b=1$, then the ratio of the two sizes is $a=1$. We prove this fact by a contradiction argument. And by the symmetry of the configuration, we assume that $0<a<1$.

Since $b=1$, employing Lemma \ref{lem3.2}, by eliminating $z$, one deduces that
\begin{eqnarray}\label{4.5}
x-ay=\frac{1}{a^{3}}x-\frac{1}{a}y.
\end{eqnarray}

By the definitions of $y$ and $z$ in (\ref{3.1}), we see that $z>y$. By \textbf{Step 1}, we have $y>0$. From (\ref{3.2}), and $0<a<1$, then $x-ay>x-z>0$. Moreover, by
\begin{eqnarray*}
|1-\rho_{k}|^{3}=|1-\rho_{N-k}|^{3},\,\,k=1,2,\ldots,N-1,
\end{eqnarray*}
and
\begin{eqnarray*}
Im(1-\rho_{k})=-\sin(\frac{2k\pi}{N})=Im(1-\rho_{N-k}),\,\,k=1,2,\ldots,N-1.
\end{eqnarray*}
Then it enables us to obtain that
$Im(x)=Im(\sum_{1\leq k\leq N-1}\[(1-\rho_{k})/(|1-\rho_{k}|^{3})\])=0$.
Hence
\begin{eqnarray*}
x&=&Re(\sum_{1\leq k\leq N-1}\[(1-\rho_{k})/(|1-\rho_{k}|^{3})\])=\frac{1-\cos(\frac{2k\pi}{N})}{|2-2\cos\frac{2k\pi}{N}|^{\frac{3}{2}}}\cr
&=&\frac{1-\cos(\frac{2k\pi}{N})}{8|\sin\frac{k\pi}{N}|^{3}}=\frac{1}{4}\sum_{1\leq k\leq N-1}\csc(\frac{ k\pi}{N})>0, \,\,\,\,\,\,k=1,2,\ldots,N.
\end{eqnarray*}
Therefore, if $0<a<1$, then
\begin{eqnarray*}
x-ay<\frac{1}{a^{2}}(x-ay)=\frac{1}{a^{2}}x-\frac{1}{a}y<\frac{1}{a^{3}}x-\frac{1}{a}y,
\end{eqnarray*}
which contradicts with (\ref{4.5}). Hence $a=1$.

\textbf{Step 3.} We claim that if $0<b<1$, then $0<a<1$.

If this statement is false, then $a\geq1$. In view of Remark \ref{rem3.1}, if $a=1$, then $b=1$. Hence, $a>1$. Thanks to (\ref{3.2}) and (\ref{3.3}), we deduce that
\begin{eqnarray*}
x-z=\frac{b^{2}}{a}x-b^{2}a^{2}z<\frac{b^{2}}{a}x-\frac{b^{2}}{a}z=\frac{b^{2}}{a}(x-z).
\end{eqnarray*}
By \textbf{Step 1}, we have $y=x-z>0$. So $b^{2}/a>1$. Thus $b^{2}>a$, which contradicts with $a>1$ and $0<b<1$. Hence $a>1$ is impossible, which implies that $0<a<1$.

\textbf{Step 4.} We claim that if $0<b<1$, then $a> b^{2}$.

From \textbf{Step 3}, we have $0<a<1$. Employing Lemma \ref{lem3.2}, by eliminating $y$, one computes that
\begin{eqnarray}\label{4.6}
x-z=\frac{b^{2}}{a}x-b^{2}a^{2}z.
\end{eqnarray}
Combining $0<a<1$ and $z>0$, we have
\begin{eqnarray}\label{4.7}
x-z>\frac{b^{2}}{a}x-\frac{b^{2}}{a}z=\frac{b^{2}}{a}(x-z).
\end{eqnarray}
Since $y=x-z>0$, by (\ref{4.7}), if $0<b<1$, then $a> b^{2}$.

\textbf{Step 5.} Employing \textbf{Step 4}, we prove that $a>b^{1/2}$.

By (\ref{4.6}), one computes that
\begin{eqnarray*}
\frac{x}{z}=\frac{a-b^{2}a^{3}}{a-b^{2}}.
\end{eqnarray*}
Set $x=(a-b^{2}a^{3})t$, $z=(a-b^{2})t$.  By \textbf{Step 4}, then $a>b^2$ and so $t>0$. By (\ref{3.2}), and the fact $y<z$, we obtain
\begin{eqnarray*}
x-z=b^{2}(1-a^{3})t=aby<abz=ab(a-b^{2})t.
\end{eqnarray*}
which implies that
\begin{eqnarray*}
b+b^{2}a<a^{2}+ba^{3}.
\end{eqnarray*}
Thus $a>b^{1/2}$.

\textbf{Step 6.} We prove that the sizes of the two regular $N$-polygons also must be equal.

By \textbf{Step 2} and $b\leq1$, in the following, it suffices to consider the case of $0<b<1$. In fact, by (\ref{3.5}) and (\ref{3.6}) in Lemma \ref{lem3.4}, and the notation (\ref{3.1}), we have
\begin{eqnarray*}
y=\frac{1}{1+b}(x+by-abz),
\end{eqnarray*}
which implies
\begin{eqnarray}\label{4.8}
y=x-abz.
\end{eqnarray}
From (\ref{3.3}) and (\ref{4.8}), we see that $b/a^{2}=1$, which contradicts with $a>b^{1/2}$ in \textbf{Step 5}. Thus $0<b<1$ is impossible, which means $b=1$.
Then by \textbf{Step 2}, we arrive at the conclusion $a=1$, which implies that the sizes of the two regular $N$-polygons also must be equal.  $\Box$

\vskip3mm\noindent
\textbf{Acknowledgements}
\vskip2mm\noindent
Liang Ding is partially supported by research funding project of Guizhou Minzu University (GZMU[2019]QN04), Science and Technology Foundation of Guizhou Province (J[2015]2074).
Jinlong Wei is partially supported by NSF of China (11501577).

\end{document}